\documentclass{amsart}
\usepackage{amssymb,latexsym,amsmath,amsfonts,hhline}
\usepackage{pdfsync}

\newcommand{\cal}{\mathcal}

\makeatletter
\renewcommand{\subsection}{\@startsection{subsection}{2}{0mm}{-2mm}{-2mm}{\bf\normalsize}}

\makeatother

\newtheorem{formula}{}[section]
\newtheorem{definition}[formula]{Definition}
\newtheorem{corollary}[formula]{Corollary}
\newtheorem{remark}[formula]{Remark}
\newtheorem{lemma}[formula]{Lemma}
\newtheorem{theorem}[formula]{Theorem}

\def\thrm{\begin{theorem}}
\def\thrml#1{\begin{theorem}\label{#1}}
\def\ethrm{\end{theorem}}
\def\rmrk{\begin{remark}}
\def\rmrkl#1{\begin{remark}\label{#1}}
\def\ermrk{\end{remark}}
\def\dfntn{\begin{definition}}
\def\dfntnl#1{\begin{definition}\label{#1}}
\def\edfntn{\end{definition}}
\def\nmrt{\begin{enumerate}}
\def\enmrt{\end{enumerate}}
\def\tm#1{\item[{\rm (#1)}]}
\def\qtn{\begin{equation}}
\def\qtnl#1{\begin{equation}\label{#1}}
\def\eqtn{\end{equation}}
\def\lmm{\begin{lemma}}
\def\lmml#1{\begin{lemma}\label{#1}}
\def\elmm{\end{lemma}}
\def\crllr{\begin{corollary}}
\def\crllrl#1{\begin{corollary}\label{#1}}
\def\ecrllr{\end{corollary}}
\def\css{\begin{cases}}
\def\ecss{\end{cases}}

\def\proof{\noindent{\bf Proof}.\ }

\def\cA{{\cal A}}
\def\cB{{\cal B}}
\def\cC{{\cal C}}

\def\cR{{\cal R}}
\def\cS{{\cal S}}

\def\cX{{\cal X}}
\def\cY{{\cal Y}}
\def\cZ{{\cal Z}}

\def\mZ{{\mathbb Z}}

\DeclareMathOperator{\aut}{Aut}

\DeclareMathOperator{\dev}{dev}

\DeclareMathOperator{\GL}{GL}

\DeclareMathOperator{\Inn}{Inn}
\DeclareMathOperator{\inv}{Inv}

\DeclareMathOperator{\mat}{Mat}
\DeclareMathOperator{\orb}{Orb}
\DeclareMathOperator{\out}{Out}

\DeclareMathOperator{\PSL}{PSL}

\DeclareMathOperator{\rk}{rk}
\DeclareMathOperator{\SL}{SL}

\DeclareMathOperator{\Span}{Span}
\DeclareMathOperator{\sz}{Sz}

\DeclareMathOperator{\sym}{Sym}

\def\bull{\hfill\vrule height .9ex width .8ex depth -.1ex\medskip}

\def\qaq{\quad\text{and}\quad}

\newcommand{\grp}[1]{\langle {#1}\rangle}
\newcommand{\und}[1]{{\underline #1}}
\def\SG#1{{\tt SmallGroup(#1)}}

\begin{document}
\title{On non-abelian Schur groups}
\author{Ilya Ponomarenko}
\address{Steklov Institute of Mathematics at St. Petersburg, Russia}
\email{inp@pdmi.ras.ru}
\author{Andrey Vasil'ev}
\thanks{The work was partially supported by RFBR, research projects No. 14-01-00156 (first author) and No. 13-01-00505 (second author)}
\address{Sobolev Institute of Mathematics and Novosibirsk State University, Novosibirsk, Russia}
\email{vasand@math.nsc.ru}
\date{05.07.14}
\maketitle

\begin{abstract}
A finite group $G$ is called Schur, if every Schur ring over $G$ is
associated in a natural way with a regular subgroup of $\sym(G)$ that is isomorphic to~$G$. We prove that any nonabelian Schur group
$G$ is metabelian and the number of distinct prime divisors of the order of $G$ does not exceed~7.
\end{abstract}

\section{Introduction}

Let $G$ be a finite group. In what follows we are interested in a  permutation group  $\Gamma\le\sym(G)$ that contains
a subgroup $G_{right}$ induced by the right multiplications in~$G$.  Denote by $\cS$ the set of orbits of the stabilizer
of the identity $e=e_G$ in $\Gamma$. Then  one can form  a $\mZ$-submodule $\cA=\cA(\Gamma,G)$ of the group ring $\mZ G$:
\qtnl{151012a}
\cA=\Span\{\und{X}:\ X\in\cS\}
\eqtn
where $\und{X}$ is the sum of elements of $X$ in the group ring. It was observed by Schur in~\cite{Sch1933} that this
submodule  is a subring of the latter ring. The idea of Schur was to study the group $\Gamma$ by means of  the ring $\cA$ rather
than by means of the characters of~$G$. In this way he was able to prove a generalization of the Burnside theorem on a permutation group
with a regular cyclic $p$-subgroup (see~\cite[Theorem~25.3]{Wie1964}).\medskip

Following Wielandt a subring $\cA$ of the ring~$\mZ G$ is called a {\it Schur ring} or {\it S-ring} over the group $G$
if there exists  a partition $\cS$ of this group such that equality~\eqref{151012a} holds and
\qtnl{171012a}
\{e\}\in\cS\qaq \cS^{-1}=\cS
\eqtn
where $\cS^{-1}=\{X^{-1}:\ X\in\cS\}$ with $X^{-1}=\{x^{-1}:\ x\in G\}$.\footnote{One can see that the partition $\cS$ is uniquely
determined.} It is easily seen that the ring $\cA(\Gamma,G)$ defined in the above paragraph is an S-ring over~$G$. In particular, so is
the group ring  itself: in this case $\Gamma=G_{right}$.\medskip

It should be mentioned that as Wielandt wrote in~\cite[p.54]{Wie1969}: {\it
Schur had conjectured for a long time that every S-ring  $\cA$ over a group $G$ is determined by a suitable permutation
group~$\Gamma$};  in our terms this means that $\cA=\cA(\Gamma,G)$.  However, this proved to be false and first counterexamples
were found by Wielandt in~\cite[Theorem~25.7]{Wie1964}. In honour of this Schur's fallacy the S-rings arising from
permutation groups were called {\it schurian} in~\cite{Poe1974}. It was also proved there that any S-ring over a cyclic
$p$-group, $p>3$, is Schurian. Thus such a group is an example of a {\it Schur group}  in the sense of the following definition given there.

\dfntnl{060212a}
A finite group $G$ is called a Schur group, if every S-ring over~$G$ is schurian.
\edfntn

The problem of determining all Schur groups was suggested by R.P\"oschel who proved that a $p$-group, $p>3$, is Schur if and
only if it is cyclic~\cite{Poe1974}. Only about thirty years after, all cyclic Schur groups were classified in~\cite{EKP2011}. The techique
developed there were used in~\cite{EKP2012} where it was proved that any abelian Schur group belongs to one of the
several explicitly given families. It should be mentioned that in the proof of this result a knowledge of Schur groups
of small order is used. In fact, all S-rings over an arbitrary group of order $\le 47$ were enumerated by means of computer
computations. It turned out that there are non-abelian Schur groups. However, except for the above result of P\"oschel, no
general results on these groups  was known. One of the main results of this paper gives a result of this type.

\thrml{100912a}
A simple group is Schur if and only if it is abelian. In particular, any Schur group is solvable.
\ethrm

The proof of Theorem~\ref{100912a} is based on the Thompson theorem on minimal simple groups, a detailed analyses of
groups of small order (Section~\ref{241012a}),  dihedral groups (Section~\ref{241012b}) and the groups $\PSL_2(q)$
(Lemma~\ref{060212e}).  Based on Theorem~\ref{100912a} and
the above analyses, we apply standard group-theoretical
arguments to obtain the second main result of the paper.

\thrml{021012a}
Any Schur group $G$ is metabelian. Moreover,
\nmrt
\tm{1} $|\pi(G)|\le 7$ with equality only if $|G|$ is even product of seven primes,
\tm{2} if a Hall $\{2,3\}$-subgroup of $G$ is cyclic, then $G$ is metacyclic.
\enmrt
\ethrm

Concerning permutation groups we refer to~\cite{Wie1964}. For the reader convenience we collect the
basic facts on S-rings and Cayley schemes in Section~\ref{241012c}. \medskip

{\bf Notation.}  For a positive integer $n$ we denote by $C_n$, $D_n$ and $E_n$ the cyclic, dyhedral
and elementary abelian groups of order $n$, and by $S_n$ and $A_n$ the symmetric and alternating groups
on $n$ letters respectively.\footnote{The notation $E_n$ has sense only if $n$ is a prime power.}

The quaternion group and generalized quaternion groups are denoted by $Q_8$ and $Q_{2^k}$, $k\geq4$, respectively.

For $k\geq4$ we set $SD_{2^k}=\grp{a,b~\mid~a^{2^{k-1}}=b^2=1, a^b=a^{-1+2^{k-2}}}$ to be
the semidihedral group.\footnote{Sometimes this group is denoted by $QD_{2^k}$ and called the quasidihedral
group.}

For a prime $p$ and $k\geq3$ we set $M_{p^k}=\grp{a,b~\mid~a^{p^{k-1}}=b^p=1, a^b=a^{1+p^{k-2}}}$, where
$k>3$ if $p=2$.

We set $G_{16}=\grp{a,b,c~\mid~a^4=b^2=c^2=[a,b]=[a,c]=1, [b,c]=a^2}.$\footnote{This
group is the \SG{16,13} in GAP notation.}

For a group $G$ we set $e_G$ to be the identity of $G$.

For groups $H$ and $G$ we write $H\prec G$ if $H$ is isomorphic to a
section of~$G$, that is a group $A/B$ with $B\triangleleft A\le G$.

The set of prime divisors of a group $G$ is denoted by $\pi(G)$.

The group ring of $G$ over integers is denoted by $\mZ G$.

The sum of elements of a set $X\subset G$ in $\mZ G$ is denoted by $\und{X}$.

The symmetric group on a set $\Omega$ is denoted by $\sym(\Omega)$.

The regular permutation group which is induced by right (resp. left) multiplications in $G$ is denoted by $G_{right}$ (resp. $G_{left}$).

\section{Schur rings and Cayley schemes}\label{241012c}
In this section we recall some facts on association schemes and Cayley schemes that will be used throughout the paper to
prove  non-schurity of  S-rings. More details can be found in~\cite{EP09}.\medskip

Let $\Gamma$ be a transitive permutation group on a finite set~$\Omega$. Denote by $\cR$ the set of orbits in its  coordinatetwise action
on the Cartesian square~$\Omega\times\Omega$. Then one can form a $\mZ$-submodule $\cZ=\cZ(\Gamma,\Omega)$ in the ring
$\mat_{\Omega}(\mZ)$ of all integer matrices the rows and columns of which are indexed by the elements of~$\Omega$:
\qtnl{201012a}
\cZ=\Span\{A(R):\ R\in\cR\}
\eqtn
where $A(R)$ is the adjacency matrix of the binary relation $R\subset\Omega\times\Omega$. It is easily seen that $\cZ$ is a matrix ring,
the centralizer ring of the permutation group~$\Gamma$. In the special case when $\Omega=G$ is a group and $G_{right}\le\Gamma$ the
injection
\qtnl{211012a}
\rho:G_{left}\to\mat_{\Omega}(\mZ),\ g\mapsto P_g,
\eqtn
where $P_g$ is the permutation matrix associated with $g$, induces a ring isomorphism from $\cA(\Gamma,G)$
onto~$\cZ(\Gamma,\Omega)$ \cite[Theorem~28.8]{Wie1964}.\medskip

A partition $\cR$ of the set $\Omega\times\Omega$ forms an {\it association scheme} $\cX=(\Omega,\cR)$ if the module $\cZ$ defined
by~\eqref{201012a} is a ring (the {\it adjacency ring} of~$\cX$) and
\qtnl{221012a}
1_\Omega\in\cR\qaq \cR^*=\cR
\eqtn
where $1_\Omega$ is the diagonal of $\Omega\times\Omega$ and $\cR^*$ consists of all relations $R^*$ obtained
from $R\in\cR$ by interchanging of coordinates.\footnote{This definition of an association scheme
is obviously equivalent to the standard one.} In the special case when $\cZ=\cZ(\Gamma,\Omega)$ for a permutation group $\Gamma$, the
corresponding association scheme is called {\it schurian} and denoted by $\inv(\Gamma)=\inv(\Gamma,\Omega)$. For an arbitrary
association scheme we define the {\it automorphism group} by
$$
\aut(\cX)=\{\gamma\in\sym(\Omega):\ R^\gamma=R,\ R\in\cR\}
$$
where $R^\gamma=\{(\alpha^\gamma,\beta^\gamma):\ (\alpha,\beta)\in R\}$. The mappings $\inv$ and $\aut$ define a Galois correspondence between
the association schemes with transitive automorphism groups and transitive
permutation groups on $\Omega$:
$$
\cX\le\inv(\Gamma)\quad\Leftrightarrow\quad\Gamma\le\aut(\cX)
$$
where the partial order on the association schemes is induced by the natural partial order of their adjacency rings.  The schurian schemes
are closed with respect to this correspondence; in the other words, a scheme $\cX$ is schurian if and only if $\cX=\inv(\aut(\cX))$.
This observation gives a simple sufficient condition for a scheme to be non-schurian.

\lmml{241012w}
A scheme $\cX$ is non-schurian whenever  $\cX\ne\inv(\aut(\cX))$.\bull
\elmm

A special class of association schemes $\cX$ arises when $\Omega=G$ is a group and
$$
G_{right}\le\aut(\cX).
$$
In this case $\cX$ is called a {\it Cayley scheme} over the group~$G$; any relation $R\in\cR$ is the arc set of the Cayley graph  on~$G$
associated with a set  $X=\{g\in G:\ (g,e)\in R\}$. The set of all such $X$'s forms a partition $\cS$ of the group~$G$ for which the
$\mZ$-module~\eqref{151012a} is an S-ring; it is denoted by~$\cA(\cX)$. Conversely, by the definition of mapping \eqref{211012a} for any
set $X\subset G$ the matrix
$$
\rho(\und{X})=\sum_{g\in X}P_g
$$
is a \{0,1\}-matrix in $\mat_{\Omega}(\mZ)$. Therefore it is the adjacency matrix of a binary relation $R\subset G\times G$ which is invariant
with respect to $G_{right}$. One can also see that all these relations form a partition $\cR$ of the set $G\times G$ for which
conditions~\eqref{221012a} are satisfied. Thus the pair $\cX(\cA)=(G,\cR)$
is a Cayley scheme over the group~$G$.

\thrml{100909a}
The correspondence $\cA\mapsto\cX(\cA)$, $\cX\mapsto\cA(\cX)$ induces a bijection between the S-rings and Cayley schemes over
the group~$G$ that preserves the natural partial orders on these sets. Moreover, under this correspondence the S-ring $\cA$ is schurian
if and only if so is the Cayley scheme $\cX(\cA)$.\bull
\ethrm

We conclude the section by two results giving necessary conditions for a group $G$ to be Schur. Below given a group $H\le G$
we denote by $\Gamma(H,G)$ the (transitive) permutation group which is induced by the action of the group~$G$ on the set of
right $H$-cosets by the right multiplications. In the proof of the following lemma we will freely use well-known properties of
the wreath product of association schemes; all of them can be found, e.g. in~\cite{E04}.

\lmml{080212f}
Let $G$ be a Schur group, $H\le G$ and $\Gamma=\Gamma(H,G)$. Then
any fusion\footnote{A scheme $\cY$ is called a {\it fusion} of a scheme $\cX$, if $\cY\le\cX$} of the scheme $\inv(\Gamma)$ is schurian.
\elmm
\proof It is easily seen that the set of right $H$-cosets forms an imprimitivity system of the group~$G_{right}$. Therefore without
loss of generality we can assume that it is a subgroup of the wreath product\  $\Delta=H_{right}\wr\Gamma$ in imprimitive action.
This implies that
$$
G_{right}\le \Delta\le \aut(\inv(\Delta))=\aut(\cX_H\wr\cX_\Gamma)
$$
where $\cX_H=\inv(H_{right})$ and $\cX_\Gamma=\inv(\Gamma)$. It follows that $\cX_H\wr\cX_\Gamma$ is a Cayley scheme over
the group~$G$. Now, suppose on the contrary that there exists a non-schurian fusion $\cY$ of the scheme $\cX_\Gamma$. Then
$$
\cX_H\wr\cY\le \cX_H\wr\cX_\Gamma,
$$
and hence $\cX_H\wr\cY$ is also a Cayley scheme over~$G$. Moreover, this scheme is schurian if and only if so is the
scheme $\cY$ (see e.g.~\cite{E04}).
By Theorem~\ref{100909a} this implies that the S-ring $\cA(\cX_H\wr\cY)$ is not schurian. Therefore the group $G$ is not Schur.
Contradiction.\bull

When the group $H$ in Lemma~\ref{080212f} is normal in $G$, the scheme $\inv(\Gamma)$ is obviously a Cayley scheme over the
group $G/H$. Therefore our result implies that the quotient of a Schur group is also Schur. This provides a proof of one part of the
following theorem proved in~\cite{Poe1974}.

\thrml{060212b}
Any section of a Schur group is a Schur group.\bull
\ethrm

\section{Abelian and small non-abelian Schur groups}\label{241012a}

All S-rings over any group of order $n\le 31$ were enumerated by computer in~\cite{F98}; the same
result for $n\le 47$ was later announced in~\cite{PR}. As a consequence of these results one can
get Lemma~\ref{241012e} below. An independent proof of this lemma we give here, is obtained by using
the enumeration of all association schemes (and hence the Cayley schemes) of small order in~\cite{HM}.

\lmml{241012e}
A nonabelian group of order $\le 27$ is not Schur unless it is a dihedral group or is isomorphic
to $Q_8$, $A_4$, $C_3:C_4$, or $G_{16}$.
\elmm
\proof By means of computer system GAP \cite{gap} one  can find that there are exactly~$22$
nonabelian groups that are neither dihedral nor isomorphic to $Q_8$, $A_4$, $C_3:C_4$, or $G_{16}$.
All that groups are listed in the first column of the Table~\ref{tbl3}; for instance, the third
row of the left-hand side correspond to the group $Q_{16}$, the \SG{16,9} in GAP.\medskip

For each group $G$ in this list we choose an association scheme $\cX=(\Omega,\cR)$
belonging to the Hanaki-Miamoto database of small association schemes~\cite{HM} where
$\Omega=\{1,\ldots,n\}$ with $n=|G|$. The second and
third columns of the table contain the index of the scheme $\cX$ in the database and the
number $\rk(\cX)=|\cR|$ respectively.\medskip

By means of the Hanaki GAP-package for association schemes~\cite{H12}
one can find the group $\Gamma=\aut(\cX)$; the number $\rk(\Gamma)=\rk(\inv(\Gamma))$ is given
in the fourth column of the Table~\ref{tbl3}. Inspecting the table shows that $\rk(\cX)\ne\rk(\Gamma)$,
and hence
$$
\cX\ne\inv(\aut(\cX)).
$$
By Lemma~\ref{241012w} this implies that the scheme $\cX$ is non-schurian. On the other hand, a
straightforward computation in GAP shows that the group $\Gamma$ contains a
regular subgroup isomorphic to~$G$. Therefore the scheme~$\cX$ is isomorphic to
a Cayley scheme over this group. Thus the group $G$ is not Schur by Theorem~\ref{100909a}.\bull

\begin{table}[bvt]
\caption{Non-schurian schemes for small groups}\label{tbl3}
\begin{center}
\begin{tabular}[tc]{|c|c|c|c|c|c|c|c|c|}
\hhline{|----|~|----|}
$G$ & $\cX$ & $\rk(\cX)$ & $\rk(\Gamma)$ & \quad\quad\quad &
$G$ & $\cX$ & $\rk(\cX)$ & $\rk(\Gamma)$ \\
\hhline{|----|~|----|}
$[16,3]$      & $[16,59]$      & $6$        & $7$          & &
$[16,4]$      & $[16,94]$      & $7$        & $10$          \\
\hhline{|----|~|----|}
$[16,6]$      & $[16,6]$        & $3$        & $4$          & &
$[16,8]$      & $[16,6]$        & $3$        & $4$            \\
\hhline{|----|~|----|}
$[16,9]$      & $[16,59]$      & $6$        & $7$          & &
$[16,11]$      & $[16,6]$        & $3$        & $4$            \\
\hhline{|----|~|----|}
$[16,12]$      & $[16,59]$      & $6$        & $7$          & &
$[18,3]$      & $[18,41]$      & $6$        & $8$            \\
\hhline{|----|~|----|}
$[18,4]$      & $[18,41]$      & $6$        & $8$          & &
$[24,1]$      & $[24,191]$      & $7$        & $24$          \\
\hhline{|----|~|----|}
$[24,3]$      & $[24,308]$      & $8$        & $9$          & &
$[24,4]$      & $[24,304]$      & $8$        & $14$          \\
\hhline{|----|~|----|}
$[24,5]$      & $[24,299]$      & $8$        & $14$          & &
$[24,7]$      & $[24,304]$      & $8$        & $14$          \\
\hhline{|----|~|----|}
$[24,8]$      & $[24,299]$      & $8$        & $14$          & &
$[24,10]$      & $[24,304]$      & $8$        & $14$          \\
\hhline{|----|~|----|}
$[24,11]$      & $[24,308]$      & $8$        & $9$          & &
$[24,12]$      & $[24,17]$      & $4$        & $6$            \\
\hhline{|----|~|----|}
$[24,13]$      & $[24,106]$      & $6$        & $12$          & &
$[24,14]$      & $[24,299]$      & $8$        & $14$          \\
\hhline{|----|~|----|}
$[27,3]$      & $[27,382]$      & $4$        & $6$          & &
$[27,4]$      & $[27,382]$      & $4$        & $6$            \\
\hhline{|----|~|----|}
\end{tabular}
\end{center}
\end{table} 

The next two statements contain almost all known information on large Schur groups. The first of
them was proved in~\cite{Poe1974}.

\thrml{150912a}
A $p$-group with $p\ge 5$ is Schur if and only if it is cyclic.\bull
\ethrm

The following theorem is a consequence of the results obtained in \cite{EKP2011} for the cyclic case
and in \cite{EKP2012} for the abelian non-cyclic case.

\thrml{040211a}
Any abelian Schur group is contained in one of the following families:
\nmrt
\tm{1} $C_{p^k}$, $C_{pq^k}$, $C_{2pq^k}$, $C_{pqr}$, $C_{2pqr}$,
\tm{2} $E_{p^k}$ where either $p=2$ and $k\le 5$, or $p=3$ and $k\le 3$,
\tm{3} $C_p\times C_{p^k}$ where $p=2$ or $p=3$,
\tm{4} $C_{2p}\times C_{2^k}$, $E_4\times C_{p^k}$, $E_4\times C_{pq}$, $E_{16}\times C_p$
where $p\ne 2$,
\tm{5} $C_6\times C_{3^k}$, $E_9\times C_{2q}$, $E_9\times C_p$ where $p\ne 2$
\enmrt
with $p,q,r$ being distinct primes and $k\ge 0$ an integer.\bull
\ethrm

\rmrkl{080613a}
From~\cite{EKP2011} it follows that all the groups in cases (1) and (2) are Schur.
\ermrk

\crllrl{050113a}
Let $G$ be an abelian Schur group of order $n$ and $\Omega(n)$ the total number
of prime factors of the integer~$n$. Then
\nmrt
\tm{1} $\pi(G)\le 4$ with equality only if $2\in\pi(G)$ and $\Omega(n)=4$,
\tm{2} if $n$ is odd, then $\pi(G)\le 3$ with equality only if $\Omega(n)=3$.\bull
\enmrt
\ecrllr

In the following two lemmas below we prove that several small groups
are not Schur. In all cases we find a transitive permutation group $\Gamma\le\sym(\Omega)$ and
an association scheme $\cX=(\Omega,\cR)$ such that
\qtnl{261212a}
\cX\ne\inv(\Gamma)\qaq\aut(\cX)=\Gamma.
\eqtn
By Lemma~\ref{241012w} this scheme is non-schurian. To define the scheme $\cX$ we will specify the
nonreflexive orbits $R_1,\ldots,R_k$ in the coordinatewise action of the group~$\Gamma$ on the set
$\Omega\times\Omega$. Then we choose a nontrivial partition $\Pi$ of the set $\{1,\ldots,k\}$,
and define a nonreflexive element of~$\cR$ to be the union of all $R_i$'s belonging to
a class of~$\Pi$. In each case one can easily check by means of GAP that (a) $\cX$ is an association scheme
and (b) relations~\eqref{261212a} hold.

\lmml{070212j}
The Frobenius groups $E_8:C_7$ and $E_{16}:C_3$ are not Schur.
\elmm
\proof Let $G$ be one of the groups in the lemma statement. Suppose we are given a permutation group
$\Gamma$ and scheme $\cX$ for which relations~\eqref{261212a} hold  with $\Omega=G$. If, in addition,
$$
G_{right}\le\Gamma,
$$
then $\cX$
is isomorphic to a non-schurian Cayley scheme over~$G$. However, then the S-ring $\cA(\cX)$ is also
non-schurian by Theorem~\ref{100909a}. Thus the group $G$ can not be Schur. To complete the
proof let us construct the group $\Gamma$ and scheme $\cX$ in each case.\medskip

Let $G=E_8:C_7=$\SG{56,11}. The group \SG{672,1257} is solvable and can be written as a
product of two subgroups that are isomorphic to~$G$ and~$A_4$ respectively. This decomposition
is not unique and the rank of the permutation representation on the right cosets
of~$A_4$ can be $8$, $12$ or $20$. In the first case, there are
three conjugacy classes of subgroups isomorhic to $A_4$,
that produce such a representation. Denote by $\Gamma$ the permutation group that corresponds to one of them. Then obviously
$\Gamma\ge G_{right}$
and the above defined number $k$ is equal to~$7$. Moreover, the indices of the
relations~$R_i$'s can be chosen so that
$$
n_1=1,\ n_2=n_3=3,\ n_4=\cdots=n_7=12
$$
where $n_i=|R_i|/|\Omega|$, $R_2$ but not $R_3$ is the union of $7$'s complete graphs on $4$ vertices,
and $(R_4)^*=R_5$ and $(R_6)^*=R_7$. Then we are done with the partition
$\Pi=\{\{1,3\},\{2\},\{4,5\},\{6,7\}\}$.\medskip

Let $G=E_{16}:C_3=$\SG{48,50}. The group \SG{1152,154768} is solvable
and is isomorpic to
a semidirect product of~$G$ by~$\SL_2(3)$. Moreover, there is a unique (up to conjugacy)
subgroup isomorphic to $\SL_2(3)$ such that the permutation representation on the right cosets of
this group has rank~$9$ and contains no a suborbit of size~$3$. Denote the corresponding permutation
group by~$\Gamma_1$. Then $N_{\sym(\Omega)}(\Gamma_1)\simeq  E_4$. Over three proper subgroups
of $N_{\sym(\Omega)}(\Gamma_1)$ that properly contain $\Gamma_1$ there are two of rank~$6$.
Denote by $\Gamma$ any of them. Then obviously $\Gamma\ge G_{right}$
and the above defined number $k$ is equal to~$5$. Moreover, the indices of the
relations~$R_i$'s can be chosen so that
$$
n_1=1,\ n_2=2,\ n_3=12,\ n_4=n_5=16.
$$
Then we are done with the partition $\Pi=\{\{1\},\{2,3\},\{4\},\{5\}\}$.\bull

\lmml{070212i}
The Frobenius group $E_{16}:C_5$ and the group $A_5$ are not Schur.
\elmm
\proof For each group $G$ from the lemma statement we choose a subgroup~$H$. By Lemma~\ref{080212f}
it suffices to find a non-schurian fusion of the scheme $\cY=\inv(\Delta,\Omega)$ where $\Delta=\Gamma(H,G)$
and $\Omega$ is the set of all right $H$-cosets. For this we will explicitly define a permutation group
$\Gamma$ such that
$$
\aut(\cY)\le\Gamma,
$$
and association scheme $\cX$ for which relations~\eqref{261212a} hold.
Then as above $\cX$ is a required non-schurian fusion of~$\cY$ and we are done.\medskip

Let $G=E_{16}:C_5=$\SG{80,49}. Set $H$ to be a subgroup of $G$ that is generated by an involution (all these
subgroups belong to the same orbit of $\aut(G)$. Then the factor group $N_{\sym(\Omega)}(\Delta)/\Delta$
has a unique normal subgroup isomorphic to $C_4\times C_2$. Denote by $\Gamma$ the permutation group
such that
$$
\Delta\le\Gamma\le N_{\sym(\Omega)}(\Delta)
$$
and the factor group $\Gamma/\Delta$ coincides with that normal subgroup.
Then the above defined number $k$ is equal to~$9$ and the indices of the relations~$R_i$'s can be chosen so that
$$
n_1=n_2=n_3=1,\ n_4=n_5=2,\ n_6=\cdots=n_9=8,
$$
and we are done with the partition $\Pi$ the unique non-singleton class of which is equal to $\{4,5\}$.
\medskip

Let $G=A_5=$\SG{60,5}. Set $H$ to be the subgroup of $G$ that is generated by the involution $(1,2)(3,4)$.
Denote by $\Gamma$ the permutation group $N_{\sym(\Omega)}(\Delta)$. Then the above defined number $k$ is equal to~$6$ and the indices of the relations~$R_i$'s can be chosen so that
$$
n_1=1,\ n_2=n_3=n_4=4,\ n_5=n_6=8,
$$
and the relations $R_2$ and $R_3$ (not $R_4$) form simple connected graphs with the vertex set~$\Omega$.
In this case we are done with the partition $\Pi$ the unique non-singleton class of which is equal to
$\{2,3\}$.\bull\medskip

We complete the section by an auxiliary lemma which is a consequence of the above results, and
will be used in Section~\ref{110113x}.

\lmml{E2group}
Let $G$ be a Schur group and $P$ a Sylow $2$-subgroup of~$G$. Suppose that $P$ is proper, normal,
elementary abelian, and $C_G(P)=P$. Then $G\simeq A_4$ or $G\simeq E_{32}:C_{31}$.
\elmm
\proof According to the Schur--Zassenhaus Theorem \cite[p.221]{GorFG}, the group $P$ being
a normal Sylow subgroup has a complement in~$G$; in other words, there exists
a group $H\le G$ of order $|G:P|$ such that $G=PH$. In particular, $H$ is of odd order.
Taking into account that $C_G(P)=P$ and $P\simeq E_{2^k}$ for some integer $k\ge 1$, without loss
of generality we can assume that
\qtnl{110113d}
H\le \GL_k(2).
\eqtn
Moreover, by statement~(2) of Theorem~\ref{040211a} we have $2\leq k\leq 5$ because $P$ is proper and $G$ is Schur. When $k=2$, then obviously
$GL_2(2)\simeq S_3$ and $|H|=3$. Therefore, $G\simeq A_4$ and
and we are done. In the remaining cases take a nonidentity element $x\in H$ of prime order~$p$.
Then by~\eqref{110113d} we have
$$
p=\css
3\ \text{or}\ 7,        &\text{if $k=3$,}\\
3,5,\ \text{or}\ 7,        &\text{if $k=4$,}\\
3,5,7\ \text{or}\ 31,        &\text{if $k=5$.}\\
\ecss
$$
In the Table~\ref{tbl4} we list   all possible triples $(k,p,c)$ with $k=3,4,5$
and $c=|C_P(x)|$. For each such triple the fourth column of this table contains a group
$K$ which is isomorphic to a subgroup of $P\grp{x}$. All these groups except for one in the last
row, are non-Schur and the reason for this is given in the fifth column. To complete the proof
it suffices to note that in the last row $K=H$ is a Sylow $31$-subgroup of $\GL_5(2)$ that is
cyclic.\bull

\begin{table}[bvt]
\caption{Reasons for groups from Lemma 3.8 to be non-schurian}\label{tbl4}
\begin{center}
\begin{tabular}{|c|c|c|c|c|}
\hline
$k$ & $p$ & $c$ & $K$ &   \\
\hline
$3$      & $3$      & $2$        & $A_4\times C_2$ &  Lemma~\ref{241012e} \\
\hline
$3$      & $7$      & $1$        & $E_8:C_7$       &  Lemma~\ref{070212j} \\
\hline
$4$      & $3$      & $1$        & $E_{16}:C_3$    &  Lemma~\ref{070212j} \\
\hline
$4$      & $3$      & $4$        & $A_4\times C_2$ &  Lemma~\ref{241012e} \\
\hline
$4$      & $5$      & $1$        & $E_{16}:C_5$    &  Lemma~\ref{070212i} \\
\hline
$4$      & $7$      & $2$        & $E_8:C_7$       &  Lemma~\ref{070212j} \\
\hline
$5$      & $3$      & $8$        & $A_4\times C_2$ &  Lemma~\ref{241012e} \\
\hline
$5$      & $3$      & $2$        & $E_{16}:C_3$    &  Lemma~\ref{070212j} \\
\hline
$5$      & $5$      & $2$        & $E_{16}:C_5$    &  Lemma~\ref{070212i} \\
\hline
$5$      & $7$      & $4$        & $E_8:C_7$       &  Lemma~\ref{070212j} \\
\hline
$5$      & $31$      & $1$       & $E_{32}:C_{31}$ &   \\
\hline
\end{tabular}
\end{center}
\end{table}

\section{Non-abelian Schur $p$-groups}

In this section we are interested in non-abelian nilpotent Schur groups. Since any nilpotent group
is the direct product of $p$-groups and for $p\ge 5$ there are no non-cyclic Schur $p$-groups, it
is quite natural to begin with studying $2$- and $3$-groups.\medskip

We recall that the rank of an abstract finite group $G$ is the least positive integer $r=r(G)$ such
that every subgroup of $G$ is generated by $r$ elements.
When $G$ is a $p$-group for a prime $p$, then by the Burnside theorem on a basis of a $p$-group $r(G)=r$ if and only if
$E_{p^r}\prec G$ and $E_{p^{r+1}}\not\prec G$.

\lmml{140912a}
A non-abelian Schur $p$-group of rank $r\geq 3$ is isomorphic to the group~$G_{16}$.
\elmm

\proof Let $G$ be a non-abelian Schur $p$-group of rank $r\geq 3$. Then it has an elementary abelian section $U/L$ of order $p^r$ such that
$U\ne G$ or $L\ne 1$. If $L\ne 1$, then $U$ contains a normal subgroup $L'<L$ with $|L/L'|=p$, and we set $H=U/L'$. If $L=1$, then the normalizer
of $U$ in $G$ contains a subgroup $U'$ with $|U/U'|=p$, and we set $H=U'/L$. Thus in any case $G$
has a section $H$
such that $|H|=p^{r+1}$ and $r(H)=r$. By Theorem~\ref{060212b} the group $H$ is Schur. However,
by Theorem~\ref{040211a} any abelian Schur group of rank at least~$3$ is elementary abelian.
Thus $H$ is not abelian (otherwise $r(H)=r+1$). Finally, $p=2$ or $p=3$ by Theorem~\ref{150912a}.
Let us consider these two cases separately.\medskip

{\bf Case $p=3$.} Without loss of generality we can assume that $H$ is a non-abelian group of order~$81$
and rank~$3$. By means of GAP one can find that \SG{81,6}$=C_{27}:C_3$ is a unique non-abelian group
of order $81$ that has no section isomorphic to $C_9:C_3$ or $E_9:C_3$. However, by Lemma~\ref{241012e}
none of latter groups is Schur. Since the group $H$ is Schur, Theorem~\ref{060212b} implies that
$H\simeq  C_{27}:C_3$. But then $r(H)=2$ which contradicts the choice of~$H$.\medskip

{\bf Case $p=2$.} Without loss of generality we  can assume that $H$ is a non-abelian group of order~$64$
and rank~$5$. However, there are exactly nine such groups
and each of them has a subgroup isomorphic to $C_4\times E_4$. By
Theorem~\ref{040211a},
this group, and hence $H$ is not Schur. Contradiction. Thus $r=3$ or $r=4$.
\medskip

Let $r=4$. In this case $H$ is a non-abelian group of order~$32$ and rank~$4$.
However, there are exactly seven such groups
and each of them has a section isomorphic to $C_4\times E_4$,
$C_2\times D_8$ or $C_2\times Q_8$. All of these groups are not Schur: the former one
by Theorem~\ref{040211a} whereas the latter two  by Lemma~\ref{241012e}. It follows that the group $H$
is not Schur. Contradiction. \medskip

Let $r=3$.  In this case $H$ is a non-abelian group of order~$16$ and rank~$3$. However,
there are exactly four such groups that have a section isomorphic to $E_8$, namely,
$$
(C_4\times C_2):C_2,\quad C_2\times D_8,\quad C_2\times Q_8,\quad G_{16}.
$$
The first three groups are not Schur by Lemma~\ref{241012e} whereas the fourth one is Schur: this
can be checked by enumeration of all S-rings over $G_{16}$. To complete the proof we will verify that
any group $K$ of order $32$ such that $G_{16}\prec K$, is not Schur. For this
one can find that $K$ is isomorphic to one of $17$ groups of order~$32$, and each of these groups
contains at least one subgroup isomorphic to one of the groups below:
$$
C_4\times C_4,\quad C_4\times E_4,\quad  C_8 : C_2,\quad SD_{16},\quad C_2\times D_8,\quad C_2\times Q_8.
$$
All these groups are not Schur: the first three by Theorem~\ref{040211a} whereas the other four
by Lemma~\ref{241012e}. Thus the group $K$ is not Schur.\bull

\thrml{100912b}
A non-abelian $p$-group is not Schur unless $p=2$ or $p=3$, and it is isomorphic to
one of the groups below:
\nmrt
\tm{1} $Q_8$, $G_{16}$, $M_{2^k}$, $k>5$, or $D_{2^k}$, $k>2$, if $p=2;$
\tm{2} $M_{3^k}$, $k>3$, if $p=3$.
\enmrt
Morover, the groups $Q_8$, $G_{16}$, $D_{2^k}$, $2<k<6$, are Schur.
\ethrm
\proof
Let $G$ be a non-abelian $p$-group. Then $p=2$ or $p=3$ by Theorem~\ref{150912a}. Moreover,
by Lemma~\ref{140912a} we can assume that $r(G)=2$. Denote by $\Phi=\Phi(G)$ the Frattini subgroup
of the group $G$. Then $G/\Phi\simeq  E_{p^r}$ for an integer $r\ge 1$, and $r=1$ if and only if the
group $G$ is cyclic. By the above assumption this
implies that $r=2$.  We claim that the group $\Phi$ is cyclic.  Suppose on the contrary that it is
not true. Then again by the assumption we have $\Phi/\Phi_2\simeq  E_{p^2}$ where $\Phi_2=\Phi(\Phi(G))$.
Therefore $G$ contains a section $H=G/\Phi_2$ such that
\qtnl{010113a}
E_{p^2}\triangleleft H\qaq H/E_{p^2}\simeq  E_{p^2},
\eqtn
in particular, $|H|=p^4$ and $r(H)=2$. Moreover, the group $H$ is not abelian, because otherwise
$H\simeq  C_{p^2}\times C_{p^2}$, and hence $H$ is not Schur by Theorem~\ref{040211a}. Next,
there are exactly five (resp. nine) non-abelian $2$-groups (res. $3$-groups)~$H$ satisfying~\eqref{010113a}.
When $p=2$ four groups
are not Schur by Lemma~\ref{241012e} whereas the fifth group is $G_{16}$ which is impossible because
$r(G_{16})=3$. When $p=3$ each of that nine groups contains a section isomorphic to
either $E_9:C_3$ or $C_9:C_3$ which are not Schur also by Lemma~\ref{241012e}. Thus the group $H$ is
not Schur. The claim is proved.\medskip

From the above claim it follows that $G/\Phi^2$ is a group of order $p^3$ that has a quotient
isomorphic to $E_{p^2}$. Since $r=2$, it is not isomorphic to $E_{p^3}$. Moreover, by Lemma~\ref{241012e}
it is also not isomorphic to a non-abelian group of order $27$ and exponent~$3$. Therefore the exponent of
$G/\Phi^2$ equals~$p^2$.  Thus the $p$-group $G$
contains a cyclic subgroup of index~$p$. So by the Burnside theorem \cite[Theorem~1.2]{B08} one of
the following statements hold:
\nmrt
\tm{1} $p=3$ and $G\simeq  M_{3^k}$ where $k\ge 3$,
\tm{2} $p=2$ and $G$ is isomorphic to one of the groups $M_{2^k}$, $Q_{2^k}$, $SD_{2^k}$ for $k\ge 4$,
\tm{3} $p=2$ and $G\simeq  D_{2^k}$ for $k\ge 3$.
\enmrt
However, both the group $Q_{2^k}$ and the group $SD_{2^k}$ contains a subgroup isomorphic to $Q_{2^{k-1}}$
for $k\ge 5$. On the other hand, the groups $Q_{16}$ and $SD_{16}$ are not Schur by Lemma~\ref{241012e}. Thus $G$ can not be isomorphic to $SD_{2^k}$, and is isomorphic to $Q_{2^k}$ only for $k=3$.
This completes the proof of the first statement. The second follows by a computer enumeration of
all S-rings over small groups.\bull\medskip

Recall that a finite nilpotent group is a direct product of
its Sylow subgroup. Therefore, the next assertion is a direct
consequence of Theorem~\ref{040211a} on a structure of abelian
Schur groups and the last theorem.

\crllrl{nilpotent}
Let $G$ be a nilpotent Schur group. Then there is at most one prime
$p$ such that a Sylow $p$-subgroup of $G$ is non-cyclic, and
$p$ is equal to $2$ or~$3$. Moreover, $G$ is a metabelian group
with $r(G)\leqslant5$ and $|\pi(G)|\leqslant4$, and if $|\pi(G)|=4$,
then $G$ is a cyclic group of even order, which is the
product of four primes.
\ecrllr

\section{Case of dihedral group}\label{241012b}

In this section we are interested in non-Schur dihedral groups. Let us begin with a
construction producing S-rings of rank~4; in fact, this construction can be considered as
a special case of one that used in~\cite[Theorem 1.6.1]{BCN} to establish
a well-known relationship between symmetric block designs and bipartite
distance regular graphs of diameter~$3$.\medskip

Let $G$ be the {\it generalized dihedral group}
associated with abelian group~$H$, i.e $G=\grp{H,g}$ with $h^g=h^{-1}$ for all $h\in H$. Suppose that
$D\subset H$ is a {\it difference set} in~$H$; by definition this means that
\qtnl{240912a}
(\und{D^{}}\cdot\und{D^{-1}})\circ \und{A}=\lambda\und{A}
\eqtn
for a positive integer~$\lambda$ where $\circ$ is the componentwise multiplication in the group
ring $\mZ H$ and $A$ is the set of nonidentity elements of~$H$. Denote
by~$\cS$ the partition of the elements of~$G$ into four classes:
$$
\{e\},\quad A,\quad X=Dg,\quad Y=(H\setminus D)g
$$
where $e=e_G$. Then obviously $A^{-1}=A$. Taking into account that $g^2=e$, we obtain that
$X^{-1}=g^{-1}D^{-1}=Dg=X$ and similarly $Y^{-1}=Y$. Thus condition~\eqref{171012a} is
satisfied for the partition~$\cS$. Furthermore, it is easily seen that
$$
\begin{array}{rcccc}
\und{A}\cdot\und{A} &=& (n-1)e & + &(n-2)\und{A}\\
\und{A}\cdot\und{X} &=& (k-1)\und{X} & + & k\und{Y}\\
\und{A}\cdot\und{Y} &=& (n-k)\und{X} & + & (n-k-1)\und{Y}\\
\end{array}
$$
where $n=|H|$ and $k=|D|$. From~\eqref{240912a} it also follows that
$$
\und{X}\cdot\und{Y}=
\und{{Dg}}\cdot\und{{Hg\setminus Dg}}=
\und{D}\cdot\und{H}-\und{D}\cdot\und{D^{-1}}=
ke+(k+\lambda)\und{A}.
$$
Finally, $\und{U}\cdot\und{V}=(\und{V}\cdot\und{U})^{-1}=\und{V}\cdot\und{U}$ for all
$U,V\in\cS$. Thus we obtain the following statement.

\lmml{240912b}
The $\mZ G$-module $\cA=\cA(D,H)=\Span\{\und{Z}:\ Z\in\cS\}$ is an S-ring over the group~$G$.\bull
\elmm

Starting with the difference set $D$ one can construct in a standard way a symmetric $2$-design
$\cB=\dev(D)$ the points, blocks and flags of which are respectively
the elements of the group~$H$, the elements of coset $Hg$, and the pairs
$(x,y)$ with $yx^{-1}\in Dg$.\footnote{The defining property of a $2$-design
is that any pair of distinct points is incident to the same number
of blocks; the design is symmetric if the number of points is equal to
the number of blocks.} The set of flags is obviously invariant with
respect to the setwise stabilizer of the set~$H$ in the group $G_{right}$. Therefore
this stabilizer forms an automorphism group of the design~$\cB$ that
acts regularly both on the points and on the blocks.
The design is $2$-transitive (resp. flag-transitive, antiflag-transitive)
if the group $\aut(\cB)$ acts $2$-transitively on the points (resp.
acts transitively on the flags, on the anti-flags).

\lmml{280912a}
The S-ring $\cA$ is schurian if and only if the design $\cB$ is
$2$-transitive, flag-transitive and antiflag-transitive.
\elmm
\proof Suppose that the S-ring $\cA$ is schurian. Then there exists a permutation group $\Gamma\le\sym(G)$ such that
$$
G_{right}\le\Gamma\qaq A,X,Y\in\orb(\Gamma_e)
$$
where $\Gamma_e$ is the stabilizer of~$e$ in $\Gamma$.
The setwise stabilizer $\Delta$ of the set $H$ in $\Gamma$ acts in a natural way on both $H$ and $Hg$.
Moreover, it is easily seen that
$$
H_{right}\le\Delta^H \qaq A\in\orb((\Delta^H)_e).
$$
Therefore the action of the group $\Delta$ on the points of~$\cB$ is
$2$-transitive. Moreover, taking into account that $X$ and $Y$ are
orbits of $\Gamma_e$, and hence of $\Delta_e$, we conclude that
the sets
$$
\bigcup_{\delta\in\Delta}(e,X)^\delta\qaq
\bigcup_{\delta\in\Delta}(e,Y)^\delta
$$
coincide respectively with the sets of flags and anti-flags of the
design~$\cB$. Therefore~$\Delta$ is an automorphism group
of this design that acts transitively on its flags and anti-flags.
Thus the design $\cB$ is $2$-transitive, flag-transitive and
antiflag-transitive. The converse statement is proved in a similar way.\bull

Given a prime power $q=4a+3$, $a\ge 1$, there is the Paley difference set~$D$ in the additive group $H$ of
a finite field with $q$ elements, that consists of $k=(q-1)/2$ non-zero squares. In this case for
all $k\ge 5$ we obviously have
$$
1+\sqrt{k}>(q-1)/k=2.
$$
Therefore by \cite[Theorem~8.3]{Kant1969} the design $\dev(D)$ is $2$-transitive only if $q=7$ or $q=11$. So
by Lemma~\ref{280912a} the S-ring $\cA(D,H)$ is not schurian for $q\ge 13$. Since the group $G=D_{2q}$ is Schur for $q\le 11$, we obtain the following statement.

\crllrl{060212d}
Let $p\equiv 3\pmod{4}$ be a prime. Then a dihedral group of order
$2p$ is a Schur group if and only if $p\le 11$.\bull
\ecrllr

To formulate the next result we recall that a difference set is
nontrivial if neither it nor its complement is a singleton.
Two difference sets $D$ and $D'$ are isomorphic if so are the
corresponding designs $\dev(D)$ and $\dev(D')$. Given a prime
power $q$ and integer $d\ge 2$ denote by $D=S_{q,d}$ the Singer
difference set: the parameters of it are as follows
\qtnl{021012c}
(n,k,\lambda)=((q^{d+1}-1)/(q-1),(q^d-1)/(q-1),(q^{d-1}-1)/(q-1)),
\eqtn
and the points and blocks of the design $\dev(D)$ are the lines and hyperplanes of a linear space of dimension $d$ over a finite field with $q$ elements.

\crllrl{290912a}
Let $G$ be the generalized dihedral group associated with abelian group~$H$ in which there exists a nontrivial difference set. Then $G$ is a Schur group only if this difference set is isomorphic
to a Singer difference set. In particular, in this case $|H|=(q^{d+1}-1)/(q-1)$ for a prime power $q$ and integer $d\ge 2$.
\ecrllr
\proof Let $D$ be the difference set from the corollary hypothesis. Suppose that the group $G$
is  Schur. Then the S-ring $\cA(D,H)$ must be schurian. So by
Lemma~\ref{280912a} the design $\dev(D)$ is $2$-transitive and antiflag-transitive. According to \cite[Statement~35a, p.91]{D68} this implies that this
design is isomorphic to the design $\dev(S)$ where $S=S_{q,d}$ for some $d$ and $q$. Thus the difference set $D$ is isomorphic to a Singer difference
set and we are done.\bull

A lot of information on difference sets with parameters~\eqref{021012c} is contained in~\cite{JPS07}.
In particular, one can find cyclic groups of order $n=(q^{d+1}-1)/(q-1)$ where $q$ is a prime
power and $d\ge 2$, each of which admits a difference set with parameters~\eqref{021012c} that is not isomorphic
to a Singer difference set.
By Corollary~\ref{290912a} this implies that in all these cases the group $D_{2n}$ is not Schur.

\section{Proof of Theorem~\ref{100912a}}

Due to Remark~\ref{080613a}, it suffices to verify that any non-abelian simple group is not Schur. However,
such a group always contains a {\it minimal simple group}, i.e a non-abelian simple group of composite order
all of whose proper subgroups are solvable. So by Theorem~\ref{060212b} only we need is to check that any
minimal non-abelian simple group $G$ is not Schur. By the Thompson Theorem \cite[Corollary~1]{T68}
the group $G$ is isomorphic to one of the following groups:
\nmrt
\tm{T1} $\PSL_2(p)$ where $p$ is a prime, $p>3$, and $p\equiv\pm2\pmod{5}$,
\tm{T2} $\PSL_2(3^p)$ where $p$ is an odd prime,
\tm{T3} $\PSL_2(2^p)$ where $p$ is a prime,
\tm{T4} $\sz(2^p)$ where $p$ is an odd prime,
\tm{T5} $\PSL_3(3)$.
\enmrt
In each of cases~(T3)--(T5) the group $G$ contains a subgroup isomorphic to a non-Schur group $H$: the second
and third columns of Table~\ref{tblsg} below indicate the group $H$ and the reason of why it is non-Schur (the existence
of such a subgroup is clear, see, e. g., \cite{CCNPW,W09}).
\begin{table}[h]
\caption{Non-schurian subgroups of some simple groups}\label{tblsg}
\begin{center}
\begin{tabular}{|l|l|l|}
\hline
$G$                                & $H$        &  non-Schur  \\
\hline
$\PSL_2(4)$                        & $A_5$      &  Lemma~\ref{070212i} \\
\hline
$\PSL_2(8),\sz(8)$                  & $E_8:C_7$  &  Lemma~\ref{070212j} \\
\hline
$\PSL_2(32),\sz(32)$              & $D_{62}$   &  Corollary~\ref{060212d} \\
\hline
$\PSL_2(2^p),\sz(2^p)$, $p>5$    & $E_{2^p}$  &  Theorem~\ref{100912b} \\
\hline
$\PSL_3(3)$                        & $S_4$      &  Lemma~\ref{241012e} \\
\hline
\end{tabular}
\end{center}
\end{table}
Thus in all these cases the group $G$ is not schur by Theorem~\ref{060212b}. In the remaining
cases (T1) and (T2) the same conclusion immediately follows from the next lemma the proof of which is
based on a construction from~\cite{Ban1991} and the description of maximal subgroups of a
symmetric group  given in~\cite{LiebPS1987}.

\medskip

\lmml{060212e}
The group $\PSL_2(q)$ is not Schur unless $q\le 3$.
\elmm
\proof The groups $\PSL_2(2)\simeq  D_6$ and $\PSL_2(3)\simeq  A_4$ are Schur by Lemma~\ref{241012e}. If $q=2^k$ and $k$ is a composite integer
with a proper prime divisor $p$, then the group $PSL_2(2^k)$ contains a subgroup isomorphic to $PSL_2(2^p)$ arising in case (T3) of the Thompson
theorem, so the group $PSL_2(q)$ with even $q$ is not Schur unless $q=2$. Since $PSL_2(5)\simeq  A_5$ is also not Schur, in what follows we may assume
that $q$ is odd and $q>5$.\medskip

Set $G=\PSL_2(q)$. The group $\Gamma=G_{right}\Inn(G)$ acts naturally on $G$, and contains a regular subgroup $G_{right}$.
Therefore $\inv(\Gamma)$ is a Cayley scheme over~$G$ such that the basic sets of the corresponding S-ring $\cA(\Gamma,G)$
are the conjugacy classes of $G$; the partition of $G$ with these classes is denoted by~$\cC$.\medskip

According to~\cite[Sections~3,4]{Ban1991}, the Cayley scheme $\inv(\Gamma)$ contains a special fusion $\cX$
of rank~$4$. The definition of this scheme given there, implies that
the S-ring $\cA=\cA(\cX)$ satisfies to the following conditions:
\nmrt
\tm{$1^*$} $\dim_\mZ(\cA)=4$,
\tm{$2^*$} each element of $\cS=\cS(\cA)$ is the union of all classes of~$\cC$ of the same size,
\tm{$3^*$} if $q$ is odd, then $\cS=\cC$ only if $q\le 5$.
\enmrt
To prove that the group~$G$ is not Schur we will verify  that the S-ring $\cA$ is not schurian.\medskip

The group $\Delta=N_{\sym(G)}(\Gamma)$ obviously contains a regular subgroup $G_{right}$, and hence
one can form the Cayley scheme $\inv(\Delta)$ and the corresponding S-ring $\cA(\Delta,G)$. Since $\Gamma$ is normal
in $\Delta$, condition $(2^*)$ implies that the point stabilizer $\Delta_e$ leaves each class of the
partition $\cS$ fixed. Moreover,
\qtnl{180911a}
\cA(\Delta,G)\supseteq\cA.
\eqtn
Next, by the main theorem in~\cite{LiebPS1987} the group $\Delta':=G^2\cdot(\out G\times C_2)$, which appears in the diagonal case~(d)
of that paper, is maximal in $A_n\Delta'$ where $n=|G|$. 
Since obviously $\Delta'$ normalizes $\Gamma=G_{left}\times G_{right}$, it follows that $\Delta'\le\Delta$. Thus, due to~\eqref{180911a} we
have
\qtnl{080212d}
\Delta'\le\Delta\le\aut(\inv(\Delta))\le\aut(\cX).
\eqtn
However, by statement~($1^*$) the group $\aut(\cX)$ can not be $2$-transitive. So it has no a subgroup isomorphic to~$A_n$. Therefore the maximality of~$\Delta'$ together with~\eqref{080212d} imply that
$\Delta'=\Delta=\aut(\cX)$. Thus to prove that the S-ring $\cA$ is not schurian it suffices to
verify that
\qtnl{170911b}
\cA(\Delta,G)\ne\cA.
\eqtn
When $q$ is odd, this inequality immediately follows from statement~($3^*$). \bull

\section{Proof of Theorem~\ref{021012a}}\label{110113x}

Let us recall some basic facts on the Fitting subgroup $F=F(G)$ of a group~$G$. By definition
$F$ is the largest normal nilpotent subgroup of~$G$. Suppose that $G$ is solvable. Then
\qtnl{050113v}
C_G(F)\leqslant F\qaq G/F\prec\aut(F)
\eqtn
(e.g., \cite[Theorem~6.1.3]{GorFG}). Moreover, $G$ has a finite Fitting height $h=h(G)$ which means
by definition that there is a normal series
\qtnl{060113a}
1=F_0<F_1<\ldots<F_h=G
\eqtn
where $F_{i+1}/F_i$ is the Fitting subgroup of $G/F_i$; clearly, $F_1=F(G)=:F$.\medskip

Turn to the proof of Theorem~\ref{021012a}. The group $G$ being a Schur group must be solvable by
Theorem~\ref{100912a}. Moreover, by Corollary~\ref{nilpotent} we can assume that $G$ is not
nilpotent and there is at most one prime $p$ such that a Sylow $p$-subgroup of $F$ is
non-cyclic, and $p$ is equal to $2$ or~$3$. Then the group~$G$ is metabelian by
Lemmas~\ref{meta}, \ref{3group} and~\ref{2group} that will be proved later. Statement~(1)
immediately follows from Corollary~\ref{050113a}. Finally, if a Hall $\{2,3\}$-subgroup of $G$ is
cyclic, then the group $F$ and by Lemma~\ref{meta} also the group $G/F$, are direct products
of Sylow subgroups. But they are cyclic by the assumtion and Theorem~\ref{150912a}. This
proves statement~(2).\medskip

In what follows we need the following general lemma on the Fitting subgroup and
its centralizer.

\lmml{centralizers}
Let $H$ be a finite solvable group and $F:=F(H)=A\times B$ where $A$ and $B$ are characteristic
subgroups of~$H$. Then
\qtnl{060113b}
C_H(F)=C_H(A)\cap C_H(B)\qaq
H/F\leqslant H/C_H(A)F\times H/C_H(B)F,
\eqtn
in particular, $H/F$ is a subgroup of a direct product of sections $K$ and $L$ of groups $\aut(A)$
and $\aut(B)$ respectively.
\elmm
\proof The first equality in~\eqref{060113b} is obvious whereas the inclusion  holds true by the
Remak Theorem. Furthermore, $H/FC_H(A)$ is a factor group of $H/C_H(A)$ and $H/FC_H(B)$ is a
factor group of $H/C_H(B)$.\bull

In the rest of the section we always assume that $G$ is a solvable but non-nilpotent Schur
group, $F=F(G)$ and $F\ne F_2$, see~\eqref{060113a}. We will also use the following
standard notation. If $\pi$ is a set of primes, then
$O_\pi(G)$ is the greatest normal $\pi$-subgroup of $G$; for a prime $p$ we also set $O_p(G)=O_{\{p\}}(G)$
and $O_{p'}(G)=O_{\pi\setminus\{p\}}(G)$.

\lmml{meta}
Suppose that  $F$ is cyclic, Then $G/F$ is abelian.
\elmm
\proof If the group $F$ is cyclic, then $\aut(F)$ is abelian, so $G/F$ is abelian by~\eqref{050113v}.\bull

\lmml{3group}
Suppose that the Sylow $3$-subgroup $P$ of $F$ is non-cyclic. Then
$$
G\simeq E_{27}:C_{13}\quad\text{or}\quad
G\simeq E_{9}\times (C_n:C_m)
$$
where $(n,m)\in\{(p,q),(2p,q),(p,2q),(p,4)\}$, $p$ and $q$ are primes greater than $3$, and
$m$ divides $p-1$.
\elmm
\proof By the Burnside theorem on a basis of a $p$-group, the action of
$O_{3'}(G/F)$ on~$P$ induces a group isomorphic to a $3'$-subgroup $T$ of~$\GL_r(3)$ where $r=r(P)$.
On the other hand, since $F$ is not cyclic, Theorem~\ref{040211a}
implies that~$r=2$ or~$r=3$. Therefore $\pi(T)\subseteq\{2,13\}$
and $13\in\pi(T)$ only if $r=3$. We claim
\qtnl{090113a}
O_2(T)=1.
\eqtn
Indeed, otherwise $G$ contains an involution, which normalizes a subgroup
of~$P$ that is isomorphic to~$E_9$. But then $G$ contains a non-abelian
section of order~$18$, which is not dihedral. However, by
Lemma~\ref{241012e} and Theorem~\ref{060212b} this contradicts the assumption that
$G$ is a Schur group.\medskip

Without loss of generality we can assume that $G$ is not isomorphic to
$E_{27}:C_{13}$. In this case we claim that
\qtnl{090113b}
r=2\qaq T=1.
\eqtn
Indeed, suppose that this is not true. Note, that  if $T\ne 1$,  then $13\in\pi(T)$
by~\eqref{090113a}.
Thus in any case $r=3$. So by Lemma~\ref{140912a} a Sylow $3$-subgroup of~$G$ is abelian. By
Theorem~\ref{040211a}  this implies that this subgroup coincides with $E_{27}$,
$P=F=E_{27}$ and $O_3(G/F)=1$. Taking into account that $G\ne F$ we conclude that $G/F\simeq  C_{13}$.
Thus $G\simeq  E_{27}:C_{13}$. Contradiction.\medskip

Recall that $F$ is a nilpotent Schur group and $P$ is a non-cyclic Sylow $3$-subgoup. Therefore by
Corollary~\ref{nilpotent} we have
$$
F=P\times Q
$$
where $Q$ is a cyclic $3'$-subgroup of $F$. Since $r=2$, Theorem~\ref{040211a} implies that
$n:=|Q|$ belongs to the set  $\{1,2,4,p,2q\}$ where $p$ and $q$ are primes $\ge 5$. In fact,
$n>2$. Indeed, otherwise the group $\aut(Q)$ is trivial. Therefore $G/F$ is isomorphic to a
factor group of $\aut(P)$ by Lemma~\ref{centralizers}. Due to~\eqref{090113b} this implies
that $F_2/F$ is a $3$-group. Moreover, this group must be trivial, because otherwise
a Sylow $3$-subgroup of $F_2$ is normal in~$G$. Thus, $F_2=F$. Contradiction.\medskip

It follows that $n\in\{4,p,2q\}$. Moreover, by Theorem~\ref{100912b} the group
$P$ is either abelian or isomorphic to $M_{3^k}$ with $k>3$, and in the same time by
Theorem~\ref{040211a} it does not contain an abelian subgroup isomorphic to $C_3\times C_9$.
Thus
$$
P\simeq E_9.
$$
Moreover, the group $O_3(G/F)$ acts on $P$ trivially (otherwise $G$ contains a non-abelian subgroup
of order~$27$, which is not Schur by Lemma~\ref{241012e}). Together with~\eqref{090113b} this
shows that  $F_2\leq C_G(P)$. Therefore by Lemma~\ref{centralizers} we have $C_G(Q)=F$ and
$G/F=G/C_G(Q)\leq\aut(Q)$. Hence $G=F_2$, and the group $G/F$ is cyclic and acts on $P$ trivially.
It follows
$$
G\simeq E_9\times (Q:K)
$$
for some $K\le\aut(Q)$. Thus the required statement follows from Theorem~\ref{040211a}
because $E_9\times K$ should be a Schur group.\bull

\lmml{2group}
Suppose that the Sylow $2$-subgroup $P$ of $F$ is non-cyclic. Then
$G$ is metabelian. Furthermore, if $r(P)\geq 3$, then
$$
G\simeq E_{32}:C_{31}\quad\text{or}\quad G\simeq P\times(C_p:C_q)
$$
where $P$ is one of the group $E_8$, $G_{16}$, $E_{16}$, the numbers $p$ and $q$ are
primes, $p>3$, and $q$ divides $p-1$.
\elmm
\proof Suppose first that $F=P$. If $r(P)\ge 3$ and $P$ is abelian, then $P\simeq E_{2^k}$, $k\ge 3$, by Theorem~\ref{100912b}. Since the group
$G_{16}$ has no subgroup isomorphic to $E_8$, Lemma~\ref{140912a} implies that any Sylow subgroup of $G$ that contains $P$, is abelian.
Due to~\eqref{050113v} we conclude that $P$ is a Sylow $2$-subgroup of~$G$. Thus the required statement follows from Lemma~\ref{241012e}.
If $r(P)\ge 3$ and $P$ is not abelian, then Lemma~\ref{140912a} implies that $P\simeq G_{16}$ and $P$ is again a Sylow $2$-subgroup of~$G$.
However, this is impossible by Lemma~\ref{241012e} applied for $G/\Phi(P)$ and $P/\Phi(P)\simeq E_8$. Thus $r(P)=2$.
If $P\simeq E_4$ or $P\simeq Q_8$, then $\aut(P)\simeq S_3$. Due to~\eqref{050113v} this implies that $|G|\le 24$,
and the required statement follows from Lemma~\ref{241012e}. By Theorems~\ref{100912b} and~\ref{040211a} in the remaining cases $P$ is isomorphic
to one of the following groups
\qtnl{130113a}
C_2\times C_{2^k}\ (k>1),\quad  M_{2^k}\ (k>5),\quad D_{2^k}\ (k>2).
\eqtn
Moreover $\aut(P)$ is a $2$-group by \cite[Theorem~34.8]{B08} and \cite[Theorem~3.5]{Cur07}. This implies that $G$ is nilpotent in
contrast to our assumption.\medskip

Let now $P$ be a proper subgroup of $F$. Then by Corollary~\ref{nilpotent} we have
\qtnl{130113b}
F=P\times Q
\eqtn
where $Q$ is a nontrivial cyclic group of odd order. By Theorem~\ref{040211a} this implies that
$2\le r\le 4$ where $r=r(P)$. Suppose first that $r=3$ or $r=4$. Then by that theorem we can assume that
$$
P\in\{E_8,G_{16},E_{16}\}\qaq Q=C_p
$$
for some odd prime $p$. By the Burnside theorem on a basis the group $O_{2'}(G/F)$ acts on $P/\Phi(P)$.
If this action is not trivial, then by Lemma~\ref{E2group} the group $G$ has a non-abelian section
of order $24$ and rank~$3$, which is impossible by Lemma~\ref{241012e}. Thus
\qtnl{130113w}
O_{2'}(G/F)\leq C_G(P)F.
\eqtn
On the other hand, if $O_{2}(G/F)\not\leq C_G(P)F$, then ф Sylow $2$-subgroup of $G$ is non-abelian and
includes $P$ as a proper subgroup. By Lemma~\ref{140912a} this implies that this subgroup is isomorphic
to~$G_{16}$ which is impossible because the latter group does not contain a subgroup isomorphic $E_8$.
Therefore
$$
O_{2}(G/F)\leq C_G(P)F.
$$
Together with \eqref{130113w} this shows that $F_2\leq C_G(P)F$. Therefore Lemma~\ref{centralizers}
implies that $C_G(Q)=F$ and $G/F=G/C_G(Q)\leq\aut(Q)$. It follows that $G=F_2$ and that $G/F$ is a cyclic
group acting on $P$ trivially. Applying also Lemma~\ref{241012e}
and Theorem~\ref{040211a}, we obtain that $G\simeq P\times (C_p:C_q)$, where $q$ and $r$ are primes,
$p>3$ and $q$ divides $p-1$.\medskip

Let now $r=2$. Then by Theorem~\ref{040211a} the group $P$ is isomorphic to $E_4$, $Q_8$ or one of the groups
in~\eqref{130113a}, and $Q=C_p$, or $P=E_4$ and $Q=C_{pq}$, where $p$ and $q$ different primes.
Let us consider three cases depending on the group~$P$. Below we set $K=G/F$.\medskip

{\bf Case 1:} $P\not\simeq E_4$ and $P\not\simeq Q_8$. In this case $\aut(P)$ is a
$2$-group (see above) and $Q\simeq C_p$ for an odd prime~$p$. By~\eqref{060113b} this implies
that the group $C_G(Q)F/F$ is isomorphic to a section of $2$-group $G/C_G(P)F\prec\aut(P)$. So
the preimage of $C_G(Q)F/F$ in $G$ is a nilpotent normal subgroup. Therefore it is contained in~$F$,
and hence $C_G(Q)\leq F$. It follows that
$$
C_G(Q)=F\qaq K\prec\aut(C_p).
$$
So $K$ is a cyclic group of order dividing~$p-1$. This proves $G$ is metabelian, as required, whenever
$P$ is abelian. So we may assume that $P$ is isomorphic to $M_{2^k}$, $k>5$, or $D_{2^k}$, $k>2$. If
$K$ has an odd order, then
$$
C_G(P)F=G=P\times O_{2'}(G),
$$
where $O_2'(G)$ is a Frobenius group $C_p:C_q$ with $q>3$. In particular, $G$ is metabelian and we are
done. Assume that $K$ is of even order. Then a Sylow $2$-subgroup $S$ of $G$ contains a proper normal subgroup
isomorphic to $P$. Therefore
$$
S\simeq M_{2^l}\quad\text{or}\quad S\simeq D_{2^l}
$$
for some integer $l>k$. However, the first case is impossible because all proper subgroups of $M_{2^l}$ are
abelian \cite[p.~29]{B08}. In the second case, $S$ contains exactly two conjugacy classes of involutions,
so index of $P$ in $S$ should be $2$ and $l=k+1$. Denote by
$C$ the cyclic subgroup of order $2^{k-1}$ in $P$. The subgroup $C$ is characteristic in $P$, so it
is normal in $G$. Put $T=C\times Q$. Then $T$ is abelian. On the other hand, $G/T$ a central
extension of the group $F/T$ of order $2$ by a cyclic group $K$, so $G/T$ is abelian. It provides
that $G$ is metabelian.\medskip

{\bf Case 2: $P\simeq Q_8$.} By Theorem~\ref{100912b} none of possible Schur $2$-groups contains a proper
subgroup isomorphic to~$P$. Therefore the group $K$ has an odd order. Moreover, by Lemma~\ref{241012e}
the order of this group is not a multiple of $3$. Since the factor group $G/C_G(P)F$ must be isomorphic to a
section of $\aut(Q_8)\simeq S_3$, we conclude that
$$
G=C_G(P)F.
$$
By Lemma~\ref{centralizers} this shows that $C_G(Q)=F$. Thus from Theorem~\ref{040211a} it follows that
$G\simeq P\times (C_p:Z_q)$, where $p$ and $q$ are primes, $q>3$, and $q$ divides $p-1$. In particular, $G$
is metabelian.\medskip

{\bf Case 3: $P\simeq E_4$.} In this case the group $F$ is abelian, so it sufficient to prove that
so is the group~$K$. If $C_G(P)\leq F$, then $K\leq S_3$. Moreover, $K$ is not isomorphic $S_3$,
because otherwise $G/Q$ is a non-abelian and non-dihedral group of order $24$, which is impossible by Lemma~\ref{241012e}.
Since $G$ is not nilpotent, this implies that $K$ is a cyclic group of order $3$ as required. On the other
hand, if $C_G(Q)\leq F$, then $K\prec\aut(Q)$. Since $Q$ is cyclic, the group $K$ is again abelian, and we are
done. Thus we can assume that
$$
C_G(P)\not\leq F\qaq C_G(Q)\not\leq F.
$$
In this case $K$ contains nontrivial subgroups $K_1=C_G(P)/F$ and
$K_2=C_G(Q)/F$, and by Lemma~\ref{centralizers} also the direct product $M=K_1\times K_2$. Since
$$
K_1\prec G/C_G(Q)\prec\aut(Q)\qaq K_2\prec G/C_G(P)\prec S_3,
$$
the group $K_1$ is abelian whereas the group $K_2$ is cyclic of order $3$. Therefore
the group $M$ is abelian, and we may assume that $K\ne M$. Then $|K/M|=2$ because $G/C_G(P)\prec S_3$.
Moreover, the group $K_1$ is of odd order, for otherwise the group $G/Q$ contains
a non-abelian and non-dihedral group of order $24$, which is impossible by Lemma~\ref{241012e}.
Therefore, $M$ is a $2'$-group. This implies that $K$ is a semidirect product of $M$ and a
subgroup $N$ of order $2$. Moreover, since $G/C_G(Q)$ is abelian, either $K$ is abelian, or $N$
acts nontrivially on $C_G(Q)$. But in the latter case $K$ contains a subgroup $C_G(Q)N$
isomorphic to $S_3$. Its preimage in $G/Q$ is a non-abelian and non-dihedral subgroup of order 24;
a contradiction. Thus, $K$ is abelian, as required. \bull\medskip

{\bf Acknowledgment.} We would like to thank Prof. M.~Muzychuk for the fruithful discussions
from which Lemmas~\ref{060212d} and~\ref{060212e} arose. We are also thankful to Prof.
M.~Klin and Dr. S.~Reichard for providing computational data on S-rings over small groups.


\begin{thebibliography}{99}

\bibitem{Ban1991}
E.~Bannai, \emph{Subschemes of some association schemes}, J. of Algebra \textbf{144} (1991), 167--188.

\bibitem{B08}
Y.~Berkovich, \emph{Groups of prime power order. Vol. 1}, de Gruyter Expositions in Mathematics, 46, Walter de Gruyter GmbH \& Co.
KG, Berlin, 2008.

\bibitem{BCN}
A.~E.~Brouwer, A.~M.~Cohen, A.~Neumaier, {\em Distance-regular graphs}.
Ergebnisse der Mathematik und ihrer Grenzgebiete, 3. Folge, 18. Berlin etc.:
Springer-Verlag, 1989.

\bibitem{JPS07}
D.~Jungnickel, A.~Pott, and  K.~W.~Smith, \emph{Difference Sets}. In: Colbourn, C.J., Dinitz, J.H. (eds) Handbook of Combinatorial Designs,
2nd edn., pp. 419--435. Chapman, Hall/CRC, Boca Raton (2007).

\bibitem{CCNPW}
J.~H.~Conway, R.~T.~Curtis, S.~P.~Norton, R.~A.~Parker, and R.~A.~Wilson,
\emph{An ATLAS of Finite Groups}, Oxford University Press, Oxford, 1985.

\bibitem{Cur07}
M.~J.~Curran, \emph{The automorphism group of a splitmetacyclic 2-group}, Arch. Math., \textbf{89}
(2007), 10--23.

\bibitem{D68}
P.~Dembowski, \emph{Finite geometries}, Ergebnisse der Mathematik und
ihrer Grenzgebiete, Band 44 Springer-Verlag, Berlin-New York, 1968.

\bibitem{E04}
S.~Evdokimov, \emph{Schurity and separability of association schemes}, Thesis (2004), 1--155.

\bibitem{EKP2011}
S.~Evdokimov, I.~Kov\'acs and, I.~Ponomarenko, \emph{Characterization of cyclic Schur groups},
Algebra and Analysis, \textbf{25} (2013), 61--85.

\bibitem{EKP2012}
S.~Evdokimov, I.~Kov\'acs, and I.~Ponomarenko, \emph{ On schurity of finite abelian groups},
{\tt arXiv:1309.0989 [math.GR]} (2013), 1--20.

\bibitem{EP09}
S.~Evdokimov, and I.~Ponomarenko, \emph{Permutation group approach
to association schemes}, Eur. J. Comb., \textbf{30} (2009),
1456-1476.

\bibitem{F98}
F.~Fiedler, \emph{Enumeration of Cellular Algebras Applied to Graphs with Prescribed Symmetry}, Master s thesis, Technische Universit\"at Dresden, 1998.

\bibitem{gap}
\emph{The GAP Group, GAP  Groups, Algorithms, and Programming}, Version
4.4.10; 2007, http://www.gap-system.org.

\bibitem{GorFG}
D.~Gorenstein, \emph{Finite Groups}, Harper \& Row, New York, 1968.

\bibitem{HM}
A.~Hanaki, and I.~Miyamoto, {\em Classification of association schemes with
small number of vertices}, published on web (http://kissme.shinshu-u.ac.jp/as/).

\bibitem{H12}
A.~Hanaki, {\em Elementary Functions for Association Schemes on GAP},
published on web (http://kissme.shinshu-u.ac.jp/as/gap/), 2012.

\bibitem{Kant1969}
W.~M.~Kantor, \emph{2-transitive symmetric designs}, Transactions of AMS
\textbf{146} (1969), 1--28.



\bibitem{LiebPS1987}
M.~W.~Liebeck, Ch.~E.~Praeger, and J.~Saxl, \emph{A classification of the maximal
subgroups of the finite alternating and symmetric groups}, J. of Algebra
\textbf{111} (1987), 365--383.

\bibitem{PR}
Ch.~Pech, and S.~Reichard, \emph{Enumerating set-orbits}, in: M.Klin, G.A.Jones, A.Jurisic, M.Muzychuk, I.Ponomarenko (Eds.), Algorithmic Algebraic
Combinatorics and Grobner Bases, Springer, 2009, 137--150.

\bibitem{Poe1974}
R.~P\"oschel, \emph{Untersuchungen von s-ringen insbesondere im gruppenring von
$p$-gruppen}, Math. Nachr. \textbf{60} (1974), 1--27.

\bibitem{Sch1933}
I.~Schur, \emph{Z\"ur theorie der einfach transitiven permutationgruppen}, S.-B.
Preus Akad. Wiss. Phys.-Math. Kl. (1933), 598--623.


\bibitem{T68}
J.~G.~Thompson, \emph{Nonsolvable finite groups all of whose local subgroups are solvable}, Bull. Amer. Math. Soc.,
\textbf{74} (1968), 383-437.

\bibitem{W09}
R.~Wilson, \emph{The finite simple groups}, Graduate Texts in Mathematics, no.
251, Springer, 2009.

\bibitem{Wie1964}
H.~Wielandt, \emph{Finite permutation groups}, Academic press, New York - London, 1964.

\bibitem{Wie1969}
H.~Wielandt, \emph{Permutation groups through invariant relations
and invariant fun\-cti\-ons}, Lect. Notes Dept. Math. Ohio St. Univ.,
Columbus, 1969.

\end{thebibliography}
\end{document}